\documentclass[11pt, a4paper]{article}

\usepackage{epsfig}
\usepackage{amssymb}
\usepackage{amsthm}
\usepackage{amsmath}
\usepackage{float}
\usepackage{multicol}
\usepackage{multirow}
\usepackage{cite}
\usepackage{enumerate}
\usepackage{xcolor}
\usepackage{tabu}
\def\ms{\medskip}
\def\nt{\noindent}

\definecolor{vividviolet}{rgb}{0.62, 0.0, 1.0}

\def\rsq{\hspace*{\fill}$\Box$\medskip}

\def\Z{\mathbb Z}
\newtheoremstyle{de}
  {10pt}          
  {10pt}  
  {\rm}  
  {}
  {\bf}  
  {. }    
  { }    
  {}     
\theoremstyle{de}

\newtheorem{example}{Example}[section]

\newtheoremstyle{theorem}
  {10pt}          
  {10pt}  
  {\it}  
  {}
  {\bf}  
  {. }    
  { }    
  {}     
\theoremstyle{theorem}

\newtheorem{theorem}{Theorem}[section]

\newtheorem{corollary}[theorem]{Corollary}

\numberwithin{equation}{section}
\def\Z{\mathbb{Z}}

\usepackage[mathlines]{lineno}
\modulolinenumbers[2]
\newcommand*\patchAmsMathEnvironmentForLineno[1]{%
\expandafter\let\csname old#1\expandafter\endcsname\csname #1\endcsname  \expandafter\let\csname oldend#1\expandafter\endcsname\csname end#1\endcsname  \renewenvironment{#1}%
{\linenomath\csname old#1\endcsname}%
{\csname oldend#1\endcsname\endlinenomath}}%
\newcommand*\patchBothAmsMathEnvironmentsForLineno[1]{%
\patchAmsMathEnvironmentForLineno{#1}%
\patchAmsMathEnvironmentForLineno{#1*}}%
\AtBeginDocument{%
\patchBothAmsMathEnvironmentsForLineno{equation}%
\patchBothAmsMathEnvironmentsForLineno{align}%
\patchBothAmsMathEnvironmentsForLineno{flalign}%
\patchBothAmsMathEnvironmentsForLineno{alignat}%
\patchBothAmsMathEnvironmentsForLineno{gather}%
\patchBothAmsMathEnvironmentsForLineno{multline}%
}

\begin{document}
\begin{center}
{\mathversion{bold}\Large \bf On  local antimagic chromatic numbers of the join of two special families of graphs}

\bigskip
{\large  Gee-Choon Lau$^a$, Wai Chee Shiu$^b$ }\\

\medskip

\emph{{$^a$}77D, Jalan Suboh, 85000 Segamat, Johor, Malaysia}\\
\emph{geeclau@yahoo.com}\\

\medskip
\emph{{$^b$}Department of Mathematics,}\\
\emph{The Chinese University of Hong Kong,}\\
\emph{Shatin, Hong Kong, P.R. China.}\\
\emph{wcshiu@associate.hkbu.edu.hk}\\

\end{center}

\begin{abstract}
It is known that null graphs and 1-regular graphs are the only regular graphs without local antimagic chromatic number.  In this paper, we  use matrices of size $(2m+1) \times (2k+1)$ to completely determine the local antimagic chromatic number of the join of null graphs, $O_m, m\ge 1,$ and 1-regular graphs of odd components, $(2k+1)P_2$, $k\ge 1$. Consequently, we obtained infinitely many (possibly disconnected or regular) tripartite graphs with local antimagic chromatic number 3. 
\ms


\end{abstract}

\baselineskip18truept
\normalsize

\section{Introduction}
Let $G=(V, E)$ be a connected graph of order $p$ and size $q$.
A bijection $f: E\rightarrow \{1, 2, \dots, q\}$ is called a \textit{local antimagic labeling}
if $f^{+}(u)\neq f^{+}(v)$ whenever $uv\in E$,
where $f^{+}(u)=\sum_{e\in E(u)}f(e)$ and $E(u)$ is the set of edges incident to $u$.
The mapping $f^{+}$ which is also denoted by $f^+_G$ is called a \textit{vertex labeling of $G$ induced by $f$}, and the labels assigned to vertices are called \textit{induced colors} under $f$.
The \textit{color number} of a local antimagic labeling $f$ is the number of distinct induced colors under $f$, denoted by $c(f)$.  Moreover, $f$ is called a \textit{local antimagic $c(f)$-coloring} and $G$ is {\it local antimagic $c(f)$-colorable}. The \textit{local antimagic chromatic number} $\chi_{la}(G)$ is defined to be the minimum number of colors taken over all colorings of $G$ induced by local antimagic labelings of $G$~\cite{Arumugam}. Let $G+H$ and $mG$ denote the disjoint union of graphs $G$ and $H$, and $m$ copies of $G$, respectively. For integers $c < d$, let $[c,d] = \{n\in\Z\;|\; c\le n\le d\}$. Very few results on the local antimagic chromatic number of regular graphs are known (see~\cite{Arumugam, LauLiShiu}).
Throughout this paper, we let $V(aP_2\vee O_m) = \{u_i, v_i, x_j\;|\; 1\le i\le a, 1\le j\le m\}$ and $E(aP_2 \vee O_m) = \{u_ix_j, v_ix_j, u_iv_i\;|\; 1\le i\le a, 1\le j\le m\}$. We also let $V(a(P_2\vee O_m)) = \{u_i, v_i, x_{i,j}\;|\; 1\le i\le a, 1\le j\le m\}$ and $E(a(P_2\vee O_m)) = \{u_ix_{i,j}, v_ix_{i,j}, u_iv_i\;|\; 1\le i\le a, 1\le j\le m\}$.

\ms \nt In~\cite{Haslegrave}, the author proved that all connected graphs without a $P_2$ component admit a local antimagic labeling. Thus, $O_m, m\ge 1$ and $aP_2, a\ge 1$ are the only families of regular graphs without local antimagic chromatic number.  In~\cite{Arumugam}, it was shown that $\chi_{la}(aP_2\vee O_1) = 1$ for $a\ge 1$. In the following sections, we extend the ideas in~\cite{LSPN-odd} to  prove that $\chi_{la}((2k+1)P_2 \vee O_m) = 3$ for all $k\ge 1$, $m\ge 2$. Moreover, we also obtain other families of (possibly disconnected or regular) tripartite graphs with local antimagic chromatic number 3.

\section{$(4n+1)\times (2k+1)$ matrix}\label{sec:2k+1}

Consider the graph $(2k+1)(P_2 \vee O_{2n})$ of order $(2k+1)(2n+2)$ and size $(2k+1)(4n+1)$.
Thus, for $i\in [1,2k+1]$ and $j\in [2,n]$, we need the following $(4n+1)\times (2k+1)$ matrix with entries in $[1, (4n+1)(2k+1)]$ that correspond to $f(u_ix_{i,j})$, $f(u_iv_i)$, $f(v_ix_{i,j})$, $1\le i\le 2k+1$, $1\le j\le 2n$, bijectively.  For $n=1$, the required $5 \times (2k+1)$ matrix has entries $f(v_ix_{i,2n-1})$, $f(v_ix_{i,2n})$, $f(u_iv_i)$, $f(u_ix_{i,1})$ and $f(u_ix_{i,2})$, $1\le i\le 2k+1$, given by the middle five rows. In the following matrix, $2\le j\le n$.
\[\fontsize{8}{11}\selectfont
\begin{tabu}{|c|[1pt]c|c|c|c|c|c|[1pt]c|c}\hline
i & 1 & 2 & 3 & \cdots & k & k+1 & \mbox{common diff.} & \\\tabucline[1pt]{-}
\multirow{2}{1.2cm}{$f(u_ix_{i,1})$} &k+1+ & k +  & k-1+  & \multirow{2}{0.3cm}{$\cdots$} & 2 +& 1 + & \multirow{2}{0.4cm}{$-1$} & \\
 & n(8k+4)  &  n(8k+4) & n(8k+4)  &  &  n(8k+4) &  n(8k+4) &  &\\\hline
\multirow{2}{1.2cm}{$f(u_ix_{i,2})$} & -3k-1 + & -3k + & -3k+1 + & \multirow{2}{0.3cm}{$\cdots$}  & -2k-2 +& -2k-1 + & \multirow{2}{0.4cm}{$+1$}  &  \\
& n(8k+4)  &  n(8k+4) & n(8k+4)  &  &  n(8k+4) &  n(8k+4) &  &\\\hline
 \vdots & \vdots & \vdots & \vdots & \cdots & \vdots & \vdots &  &\\\hline
\multirow{2}{2cm}{$f(u_ix_{i,2n-2j+1})$} &k+1+ & k +  & k-1+  & \multirow{2}{0.3cm}{$\cdots$} & 2 +& 1 + & \multirow{2}{0.4cm}{$-1$} & \\
 & j(8k+4)  &  j(8k+4) & j(8k+4)  &  &  j(8k+4) &  j(8k+4) & &\\\hline
\multirow{2}{2cm}{$f(u_ix_{i,2n-2j+2})$} & -3k-1 + & -3k + & -3k+1 + & \multirow{2}{0.3cm}{$\cdots$}  & -2k-2 +& -2k-1 + &  \multirow{2}{0.4cm}{$+1$} &   \\
& j(8k+4)  &  j(8k+4) & j(8k+4)  &  &  j(8k+4) &  j(8k+4) & &\\\hline
\vdots & \vdots & \vdots & \vdots & \cdots & \vdots & \vdots &  \\\hline
f(u_ix_{i,2n-1}) & 10k+5 & 10k+3 & 10k+1 & \cdots & 8k+7 & 8k+5 & -2 &   \\\hline
f(u_ix_{i,2n}) &  5k+3 & 5k+4 & 5k+5 & \cdots & 6k+2 & 6k+3 & +1 & \\\tabucline[1pt]{-}
 f(u_iv_i) & 1 & 2 & 3 & \cdots & k & k+1 & +1 & \\\tabucline[1pt]{-}
 f(v_ix_{i,1}) & 3k+2 & 3k+3 & 3k+4 & \cdots & 4k+1 & 4k+2 & +1 &\\\hline
f(v_ix_{i,2}) &   8k+4 & 8k+2 & 8k & \cdots & 6k+6 & 6k+4 & -2 &\\\hline
\vdots & \vdots & \vdots & \vdots & \cdots & \vdots & \vdots &  &\\\hline
\multirow{2}{1.5cm}{$f(v_ix_{i,2j-1})$} & -5k-2 + & -5k-1 +  & -5k +& \multirow{2}{0.3cm}{$\cdots$} & -4k-3 + & -4k-2 +&  \multirow{2}{0.4cm}{$+1$} &  \\
 & j(8k+4)  &  j(8k+4) & j(8k+4)  &  &  j(8k+4) &  j(8k+4) & &\\\hline
\multirow{2}{1.2cm}{$f(v_ix_{i,2j})$} & -k + & -k-1 + & -k-2 + & \multirow{2}{0.3cm}{$\cdots$}  & -2k+1 + & -2k + & \multirow{2}{0.4cm}{$-1$} & \\
 & j(8k+4)  &  j(8k+4) & j(8k+4)  &  &  j(8k+4) &  j(8k+4) &  \\\hline
 \vdots & \vdots & \vdots & \vdots & \cdots & \vdots & \vdots &  \\\hline
\multirow{2}{1.2cm}{$f(v_ix_{i,2n-1})$} &   -5k-2 + & -5k-1 +  & -5k +& \multirow{2}{0.3cm}{$\cdots$} & -4k-3 + & -4k-2 +& \multirow{2}{0.4cm}{$+1$} &   \\
 & n(8k+4)  &  n(8k+4) & n(8k+4)  &  &  n(8k+4) &  n(8k+4) &\\\hline
 \multirow{2}{1.2cm}{$f(v_ix_{i,2n})$} & -k + & -k-1 + & -k-2 + &  \multirow{2}{0.3cm}{$\cdots$} & -2k+1 + & -2k + & \multirow{2}{0.4cm}{$-1$} & \\
 & n(8k+4)  &  n(8k+4) & n(8k+4)  &  &  n(8k+4) &  n(8k+4) &\\\hline
\end{tabu}\]
\[\fontsize{9}{11}\selectfont
\begin{tabu}{|c|[1pt]c|c|c|c|c|c|[1pt]c|}\hline
i  & k+2 & k+3 & \cdots & 2k-1 & 2k & 2k+1 & \mbox{common diff.}\\\tabucline[1pt]{-}
\multirow{2}{1.2cm}{$f(u_ix_{i,1})$} &  2k+1+ & 2k+ & \multirow{2}{0.3cm}{$\cdots$} & k+4+ & k+3+ & k+2+ & \multirow{2}{0.4cm}{$-1$} \\
& n(8k+4)  &  n(8k+4) &  & n(8k+4) &  n(8k+4) &  n(8k+4) & \\\hline
\multirow{2}{1.2cm}{$f(u_ix_{i,2})$} & -4k-1+ & -4k+ & \multirow{2}{0.3cm}{$\cdots$} & -3k-4+ & -3k-3+ & -3k-2+    & \multirow{2}{0.4cm}{$+1$} \\
 & n(8k+4)  &  n(8k+4) &  & n(8k+4) & n(8k+4) &  n(8k+4) & \\\hline
\vdots & \vdots & \vdots & \cdots & \vdots & \vdots & \vdots  & \\\hline
\multirow{2}{2cm}{$f(u_ix_{i,2n-2j+1})$} &   2k+1+ & 2k+ & \multirow{2}{0.3cm}{$\cdots$} & k+4+ & k+3+ & k+2+ &  \multirow{2}{0.4cm}{$-1$}\\
 & j(8k+4)  &  j(8k+4) &   & j(8k+4) &  j(8k+4) &  j(8k+4) & \\\hline
\multirow{2}{2cm}{$f(u_ix_{i,2n-2j+2})$} & -4k-1+ & -4k+ & \multirow{2}{0.3cm}{$\cdots$} & -3k-4+ & -3k-3+ & -3k-2+ & \multirow{2}{0.4cm}{$+1$} \\
 & j(8k+4)  &  j(8k+4) &  & j(8k+4) & j(8k+4) &  j(8k+4) & \\\hline
\vdots & \vdots & \vdots & \cdots & \vdots & \vdots & \vdots  & \\\hline
f(u_ix_{i,2n-1}) & 10k+4 & 10k+2 & \cdots & 8k+10 & 8k+8 & 8k+6 & -2 \\\hline
f(u_ix_{i,2n}) & 4k+3 & 4k+4 & \cdots & 5k & 5k+1 & 5k+2 & +1\\\tabucline[1pt]{-}
 f(u_iv_i) & k+2 & k+3 & \cdots & 2k-1 & 2k & 2k+1 & +1\\\tabucline[1pt]{-}
f(v_ix_{i,1}) &  2k+2 & 2k+3 & \cdots & 3k-1 & 3k & 3k+1& +1 \\\hline
f(v_ix_{i,2}) & 8k+3 & 8k+1 & \cdots & 6k+9 & 6k+7 & 6k+5 & -2\\\hline
\vdots & \vdots & \vdots & \cdots & \vdots & \vdots & \vdots  & \\\hline
\multirow{2}{1.5cm}{$f(v_ix_{i,2j-1})$} & -6k-2+ & -6k-1+ & \multirow{2}{0.3cm}{$\cdots$} & -5k-5+ & -5k-4 & -5k-3 +  &\multirow{2}{0.4cm}{$+1$}   \\
 & j(8k+4)  &  j(8k+4) &   &  j(8k+4) & j(8k+4) &  j(8k+4) & \\\hline
\multirow{2}{1.2cm}{$f(v_ix_{i,2j})$} & 0+ & -1+ & \multirow{2}{0.3cm}{$\cdots$}  & -k+3+ & -k+2+ & -k+1+  & \multirow{2}{0.4cm}{$-1$}   \\
 & j(8k+4)  &  j(8k+4) &   &  j(8k+4) & j(8k+4) &  j(8k+4) &\\\hline
\vdots & \vdots & \vdots & \cdots & \vdots & \vdots & \vdots  & \\\hline
\multirow{2}{1.5cm}{$f(v_ix_{i,2n-1})$} & -6k-2+ & -6k-1+ & \multirow{2}{0.3cm}{$\cdots$} & -5k-5+ & -5k-4 & -5k-3 +  &  \multirow{2}{0.4cm}{$+1$}   \\
 & n(8k+4)  &  n(8k+4) &  &  n(8k+4) & n(8k+4) &  n(8k+4) &\\\hline
\multirow{2}{1.2cm}{$f(v_ix_{i,2n})$} &  0+ & -1+ & \multirow{2}{0.3cm}{$\cdots$}  & -k+3+ & -k+2+ & -k+1+  & \multirow{2}{0.4cm}{$-1$}  \\
 & n(8k+4)  &  n(8k+4) &  &  n(8k+4) & n(8k+4) &  n(8k+4) &\\\hline
\end{tabu}
\]

\nt We now have the following observations.
\begin{enumerate}[(1)]
\item For a fixed $j\in [2, n]$,\\
$\begin{aligned} &\quad\ \{f(u_i, x_{i, 2n-2j+1}), f(u_i, x_{i, 2n-2j+2}), f(v_i, x_{i, 2j-1}), f(v_i, x_{i, 2j})\;|\; 1\le i\le 2k+1 \}\\&=[-6k-2+j(8k+4),\ 2k+1+j(8k+4)].\end{aligned}$
\\ Thus, when $j$ runs through $[2,n]$, the integers from $10k+6$ to $2k+1+n(8k+4)$ are used. From the middle 5 rows, one may see that the integers from $1$ to $10k+5$ are used.  Thus all integers in $[1, (4n+1)(2k+1)]$ are used once.
\item For each $i\in [1,2k+1]$, the sum of the first $2n+1$ entries is $f^+(u_i) = 8n^2k+ 6nk +4n^2 +4n+k+1$.
\item For each $i\in [1,2k+1]$, the sum of the last $2n+1$ entries is $f^+(v_i) = 8n^2k+ 2nk +4n^2 +2n+k+1$.
\item Suppose $n=1$. For $j=1,2$, we let $S_j =\{f^+(x_{i,j})=f(u_ix_{i,j})+f(v_ix_{i,j})\;|\; 1\le i\le 2k+1\}$. Thus, the elements of $S_j$ form an arithmetic sequence with first term $13k+7$ and last term $11k+7$ with common difference $-1$. The sum of all the elements in $S_j$ is $(2k+1)(12k+7)$.
\item Suppose $n\ge 2$. For $j \in [1,2n]$, we also let $S_j =\{f^+(x_{i,j})=f(u_ix_{i,j})+f(v_ix_{i,j})\mid 1\le i\le 2k+1\}$.  If $j=2, 2n-1$, the elements of $S_j$ form an arithmetic sequence with first term $5k+3+n(8k+4)$ and last term $3k+3+n(8k+4)$ with common difference $-1$. For each $j\in [1,2n]\setminus \{2,2n-1\}$, $S_j = \{(4k+3)+n(8k+4), \ldots, (4k+3)+n(8k+4)\}$ with multiplicity $2k+1$. Thus, the sum of all the elements in each $S_j$, $1\le j\le 2n$, is $(2k+1)[(4k+3) + n(8k+4)]$.
\item Suppose $2k+1= (2r+1)(2s+1)$, $r,s\ge 1$. Note that, each $S_j$ is an arithmetic sequence with common difference either $0$ or $-1$. If $S_j$ is an arithmetic sequence with common difference $-1$, then $S_j$ can be partitioned into $2r+1$ blocks of size $2s+1$ such that sum of all elements in each block is $(2k+1)[(4k+3)+n(8k+4)]/(2r+1) = (2s+1)[(4k+3) + n(8k+4)]$ by the existence of $(2r+1)\times (2s+1)$ magic rectangle. If $S_j$ is an arithmetic sequence with common difference 0, then the partition is obvious. This corresponds to partition the $2k+1$ columns into $2r+1$ blocks of $2s+1$ columns so that for each $j\in [1,2n]$, $\sum [f(u_ix_{i,j}) + f(v_ix_{i,j})] =  (2s+1)[(4k+3)+n(8k+4)]$ over all the $2s+1$ columns in the same block. Moreover, the partition of $S_1$ to $S_{2n}$ may be all distinct.
\end{enumerate}

\begin{theorem}\label{thm-oddP2VO2n} For $n,k\ge 1$, $\chi_{la}((2k+1)P_2 \vee O_{2n}) = 3$. \end{theorem}

\begin{proof} Note that $\chi_{la}((2k+1)P_2 \vee O_{2n})\ge \chi((2k+1)P_2 \vee O_{2n}) = 3$. Let $G = (2k+1)(P_2 \vee O_{2n})$ be defined at the beginning of this section. We now define a bijection $f: E(G) \to [1, (2k+1)(4n+1)]$ according to the table above. Clearly, for $1\le i\le 2k+1$, $f^+(u_i) =  8n^2k+ 6nk +4n^2  +4n+k+1 > f^+(v_i) = 8n^2k+ 2nk +4n^2  +2n+k+1$. Now, for each $j\in [1, 2n]$, merging the vertices $x_{i,j}, 1\le i\le 2k+1$, to form new vertex $x_j$ of degree $4k+2$, we get the graph $(2k+1)P_2 \vee O_{2n}$. From Observations (2) to (5) above, we get that $(2k+1)P_2 \vee O_{2n}$ that admits a bijective edge labeling $f$ with
\begin{enumerate}[$(a)$]
\item $f^+(x_j) = (2k+1)[(4k+3)+n(8k+4)]$,
\item $f^+(u_i) = 8n^2k+ 6nk +4n^2 +4n+k+1$ and
\item $f^+(v_i) = 8n^2k+ 2nk +4n^2 +2n+k+1$.
\end{enumerate}
\nt Now,
\begin{align*}
(a) - (b)   &=  16k^2n-8kn^2+8k^2 +10kn-4n^2 +9k+2 \\
 & =  8kn(2k - n) + 8k^2 + 2kn + 4n(2k-n) + 9k + 2 \\
 &  = (8kn+4n)(2k-n) + 8k^2 + 2kn + 9k + 2\\
 & > 0 \quad \mbox{ if $2k\ge n$.}
\end{align*}
Otherwise, if $n\ge 2k+1$ (equivalently, $-n \le -2k-1$) we get that $(a) - (b) \le -8kn - 4n + 8k^2 + 2kn + 9k+2 \le 6k(-2k-1)+4(-2k-1)+8k^2+9k+2 < 0$. Thus, $f^+(x_j)\ne f^+(u_i)$. Similarly,
\begin{align*}
(a) - (c) & =  16k^2n-8kn^2 + 8k^2 +14kn - 4n^2 +9k+2n+2 \\
&= 8kn(2k-n) + 8k^2 + 6kn + 4n(2k-n) + 9k + 2n + 2 \\
&= (8kn+4n)(2k-n) + 8k^2 + 6kn+9k+2n+2\\
& >  0 \quad\mbox{ if $2k\ge n$.}
\end{align*}
If $n = 2k+1$, we can get $(a) - (c)  = 4k^2 + 3k > 0$. Otherwise, if $n\ge 2k+2$ (equivalently, $-n \le -2k-2$) we get that $(a) - (c) \le -16kn + 8k^2 + 6kn - 8n + 9k + 2n + 2 \le 10k(-2k-2) + 8k^2 +6(-2k-2)+9k+2 < 0$. Thus, $f^+(x_j)\ne f^+(v_i)$.

\ms\nt Consequently, $f$ is a local antimagic 3-coloring and $\chi_{la}((2k+1)P_2 \vee O_{2n})\le 3$. This completes the proof.
\end{proof}

\begin{example}\label{eg-9P2VO4} Consider $n=2, k=4$. The $9\times 9$ matrix and the $9(P_2\vee O_4)$ with the defined edge labeling are given below. For each $j\in[1,4]$, merge the vertices in $\{x_{i,j}\mid 1\le i\le 9\}$ to get the vertex $x_j$ of degree 18 gives the $9P_2 \vee O_4$ with the defined labeling as required. The induced vertex labels of $u_i, v_i, x_j$ are 205, 169, 819 respectively.

\[
\begin{tabu}{|c|[1pt]c|c|c|c|c|c|c|c|c|}\hline
i & 1 & 2 & 3 & 4 & 5 & 6 & 7 & 8 & 9  \\\tabucline[1pt]{-}
f(u_ix_{i,1}) & 77 & 76 & 75 & 74 & 73 & 81 & 80 & 79 & 78  \\\hline
f(u_ix_{i,2}) & 59 & 60 & 61 & 62 & 63 & 55 & 56 & 57 & 58  \\\hline
f(u_ix_{i,3}) & 45 & 43 & 41 & 39 & 37 & 44 & 42 & 40 & 38    \\\hline
f(u_ix_{i,4})  & 23 &  24 & 25 & 26 & 27 & 19 & 20 & 21 & 22   \\\tabucline[1pt]{-}  
f(u_iv_i) & 1 & 2 & 3 & 4  & 5 & 6 & 7 & 8 & 9  \\\tabucline[1pt]{-}
f(v_ix_{i,1}) & 14 & 15 &  16 & 17 &  18  & 10 & 11 & 12 & 13 \\\hline
f(v_ix_{i,2}) & 36 & 34 & 32 & 30 & 28  & 35 & 33 & 31  & 29  \\\hline
f(v_ix_{i,3})  & 50 & 51 & 52 & 53 & 54  & 46 & 47 & 48 & 49 \\\hline
f(v_ix_{i,4})  & 68  & 67 & 66 & 65 & 64 & 72 & 71 & 70 & 69 \\\hline
\end{tabu}\]

\begin{figure}[H]
\centerline{\epsfig{file=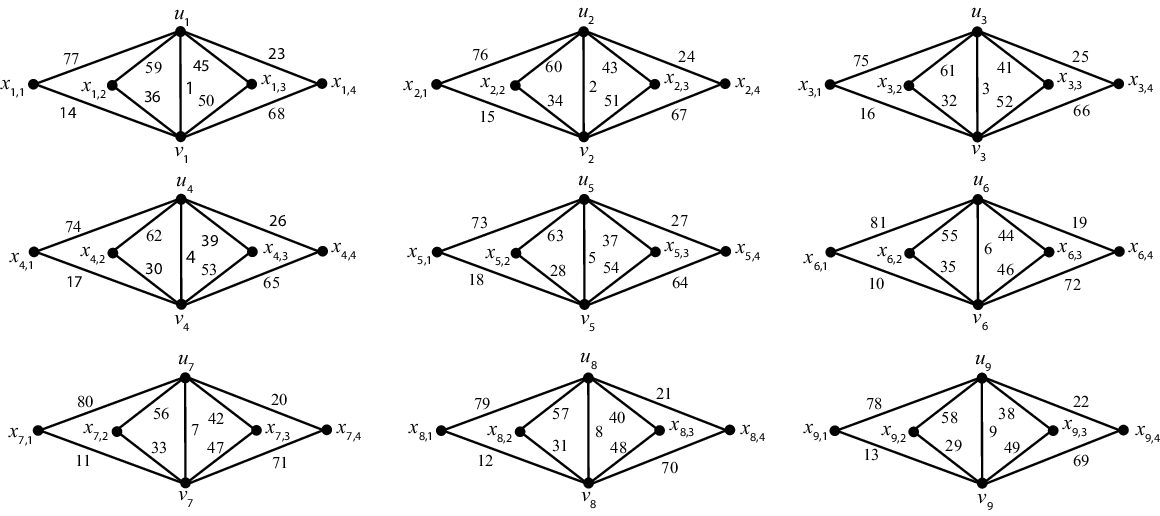, width=15.5cm}}
\caption{The graph $9(P_2\vee O_4)$ for the graph $9P_2 \vee O_4$.}\label{fig:fb9}
\end{figure}
\rsq
\end{example}

\nt Let $2k+1 = (2r+1)(2s+1)$ for $r,s\ge 1$. Consider $(2k+1)(P_2 \vee O_2)$ with the bijective edge labeling as defined  in the proof of  Theorem~\ref{thm-oddP2VO2n}. Recall that the $i$-th $P_2\vee O_2$ has two degree 3 vertices $u_i, v_i$ and two degree 2 vertices $x_{i,j}$, $j = 1,2$, with induced vertex labels $f^+(u_i) = 15k+9$, $f^+(v_i) = 11k+7$, and $f^+(x_{i,j}) = 13k+8-i$ for $1\le i\le 2k+1$.  We now have $S_j = \{f^+(x_{i,j})\mid 1\le i\le 2k+1\}$. By Observations (4) and (6) above, we can now partition $S_j$ into $2r+1$ blocks of size $2s+1$ each such that in each block, the sum of the induced vertex labels is the constant $(2s+1)(12k+7)$. Now, let $\mathcal G_2(2r+1,2s+1)$ be the set of all the graphs obtained by merging all the $2s+1$ vertices of degree 2  with induced vertex labels in the same block to get a new vertex of  degree $4s+2$ with induced vertex label $(2s+1)(12k+7)$. Since the partition of each $S_j$ may not be unique, $\mathcal G_2(2r+1,2s+1)$ may contain more than one graph (see Example~\ref{eg-G2(3,3)}).

\begin{theorem}\label{thm-G2(2r+1,2s+1)} Let $2k+1 = (2r+1)(2s+1)$ for $r,s\ge 1$. If $G\in \mathcal G_2(2r+1,2s+1)$ is defined as above, then $\chi_{la}(G) = 3$. \end{theorem}

\begin{proof} By the discussion above, we know each graph in  $\mathcal G_2(2r+1,2s+1)$ admits a local antimagic 3-coloring with distinct induced vertex labels $15k+9$, $11k+7$, $(2s+1)(12k+7)$. Thus, $\chi_{la}(G) \le 3$. Since $\chi_{la}(G) \ge \chi(G) = 3$, the theorem holds.
\end{proof}

\begin{corollary} For $r,s\ge 2$, $\chi_{la}((2r+1)[(2s+1)P_2\vee O_2]) = 3$. \end{corollary}

\begin{proof}  Let the partition of $S_1$ and $S_2$ be equal, then the resulting graph is $(2r+1)[(2s+1)P_2 \vee O_2]$.  \end{proof}

\begin{example}\label{eg-G2(3,3)} Consider $k=4$ so that $r=s=1$. We have the following $5\times 9$ matrix and a magic square $M$ of order 3.
\[
\begin{tabu}{|c|[1pt]c|c|c|c|c|c|c|c|c|}\hline
i & 1 & 2 & 3 & 4 & 5 & 6 & 7 & 8 & 9  \\\tabucline[1pt]{-}
f(u_ix_{i,1}) & 45 & 43 & 41 & 39 & 37 & 44 & 42 & 40 & 38    \\\hline
f(u_ix_{i,2})  & 23 &  24 & 25 & 26 & 27 & 19 & 20 & 21 & 22   \\\tabucline[1pt]{-}
f(u_iv_i) & 1 & 2 & 3 & 4  & 5 & 6 & 7 & 8 & 9  \\\tabucline[1pt]{-}
f(v_ix_{i,1}) & 14 & 15 &  16 & 17 &  18  & 10 & 11 & 12 & 13 \\\hline
f(v_ix_{i,2}) & 36 & 34 & 32 & 30 & 28  & 35 & 33 & 31  & 29  \\\hline
\end{tabu} \qquad
M=\begin{array}{|c|c|c|}\hline
8 & 3 & 4\\\hline
1 & 5 & 9\\\hline
6 & 7 & 2\\\hline
\end{array}
\]

\nt For $j=1,2$, partition $S_j$ into blocks $\{f^+(x_{1,j}), f^+(x_{5,j}), f^+(x_{9,j})\}$, $\{f^+(x_{2,j}), f^+(x_{6,j}), f^+(x_{7,j})\}$, $\{f^+(x_{3,j}), f^+(x_{4,j}), f^+(x_{8,j})\}$ by using the rows of $M$. We get $G=3(3P_2\vee O_2)\in \mathcal G_2(3,3)$ as required. The induced vertex labels are $69,51,165$ respectively.

\begin{figure}[H]
\centerline{\epsfig{file=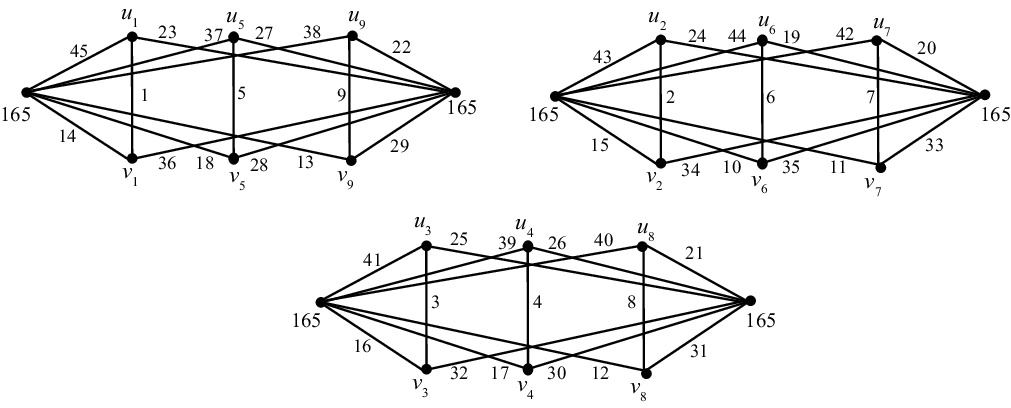, width=14cm}}
\caption{The graph $3(3P_2\vee O_2)$.}
\end{figure}

\nt If we keep the partition of $S_1$ and partition $S_2$ into blocks $\{f^+(x_{1,2}), f^+(x_{6,2}), f^+(x_{8,2})\}, \{f^+(x_{2,2}),$ $f^+(x_{4,2}), f^+(x_{9,2})\}$ and $\{f^+(x_{3,2}), f^+(x_{5,2}), f^+(x_{7,2})\}$ by using the columns of $M$, then we get a connected graph as follow with the same induced vertex labels.

\begin{figure}[H]
\centerline{\epsfig{file=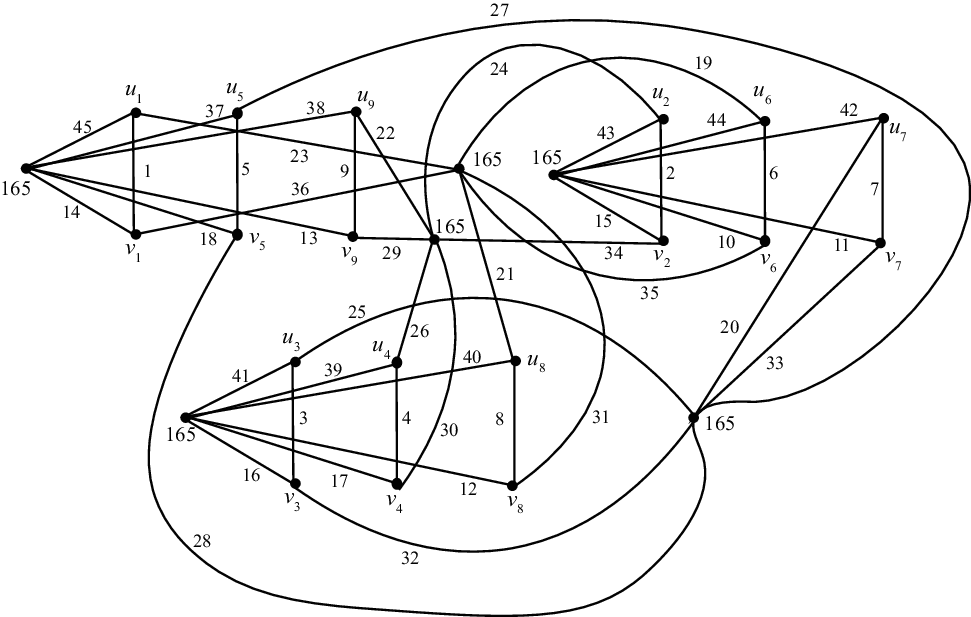, width=12cm}}
\caption{A connected graph in $\mathcal G_2(3,3)$.}
\end{figure}

\nt The degree 6 vertices are indicated with their induced vertex label 165. \rsq
\end{example}

\nt Let $2k+1 = (2r+1)(2s+1)$ for $r,s\ge 1$. Consider $(2k+1)(P_2 \vee O_{2n}), n\ge 2,$ with the bijective edge labeling as defined in the proof of Theorem~\ref{thm-oddP2VO2n}.  Similar to the discussion above for $(2k+1)(P_2 \vee O_2)$, and by Observations (5) and (6), we also have the set of graphs $\mathcal G_{2n}(2r+1,2s+1)$. Each of the graph has vertices $u_i$ and $v_i$, $1\le i\le 2k+1$, of degree $2n+1$ with induced vertex labels $f^+(u_i) = 8n^2k+6nk+4n^2+4n+k+1$, $f^+(v_i) = 8n^2k+2nk+4n^2+2n+k+1$, and $2n(2r+1)$ vertices of degree $2(2s+1)$ with induced vertex label $(2s+1)[(4k+3)+n(8k+4)]$.
Note that, since $S_1$ is a constant sequence, there are at least two different partitions of $S_1$. Thus  $\mathcal G_{2n}(2r+1,2s+1)$ contains at least two graphs.

\begin{theorem}  \label{thm-G_{2n}(2r+1,2s+1)} For $n\ge 2$, $r,s\ge 1$ and $G\in \mathcal G_{2n}(2r+1,2s+1)$, $\chi_{la}(G) = 3$. \end{theorem}

\begin{proof} By the discussion above, we know each $G\in \mathcal G_{2n}(2r+1,2s+1)$ admits a local antimagic 3-coloring with each degree $2(2s+1)$ vertex has induced vertex label $(a) = (2s+1)[(4k+3)+n(8k+4)]$ whereas for $1\le i\le 2k+1$, the vertices $u_i, v_i$ of degree $2n+1$ have induced vertex labels $(b) = f^+(u_i) =  8n^2k+6nk+4n^2+4n+k+1$ and $(c) = f^+(v_i) = 8n^2k+2nk+4n^2+2n+k+1$.

\ms\nt Clearly, $(b) > (c)$.  Now,
\begin{align*}
(a) - (b) &= 16kns-8kn^2 +2kn+8ks +8ns-4n^2+3k+6s+2\\
&= (8kn+ 4n)(2s-n) + 2kn+8ks+3k+6s+2\\
&> 0 \quad\mbox{ if $2s\ge n$.}
\end{align*}
Otherwise, if $n\ge 2s+1$ (equivalently, $-n\le -2s-1$), $(a) - (b) \le -6kn - 4n + 8ks + 3k +6s + 2\le 6k(-2s-1) + 4(-2s-1) + 8ks + 3k + 6s + 2<0$. Thus, $(a) \ne (b)$. Similarly,

\begin{align*}
(a) - (c) &= 16kns-8kn^2+6kn+8ks-4n^2 +8ns+3k+2n+6s+2\\
&= (8kn+4n)(2s-n)+6kn+8ks+3k+2n+6s+2\\
&> 0 \quad \mbox{ if $2s\ge n$.}
\end{align*}
If $n=2s+1$, $(a) - (c) =(-2k-2)(2s+1)+8ks+3k+6s+2 > 0$. Otherwise, if $n\ge 2s+2$ (equivalently, $2s-n\le -2$), $(a) -(c) \le -10kn - 6n + 8ks + 3k + 6s + 2 \le (10k+6)(-2s-2) + 8ks + 3k + 6s + 2 < 0$. Thus, $(a) \ne (c)$.

\ms\nt Consequently,  $\chi_{la}(G)\le 3$. Since $\chi_{la}(G) \ge \chi(G) = 3$, the theorem holds.
\end{proof}

\begin{corollary} For $n\ge 2, r,s\ge 1$, $\chi_{la}((2r+1)[(2s+1)P_2 \vee O_{2n}]) = 3$. \end{corollary}

\begin{proof} By using the rows of a $(2r+1)\times (2s+1)$ magic rectangle to partition $S_2$. Partition other $S_j$ so that they have the same partition as $S_2$. The resulting graph is $(2r+1)[(2s+1)P_2 \vee O_{2n}]$. \end{proof}

\begin{example}\label{eg-G4(33)} Consider $n=2, k= 4$ as in Example~\ref{eg-9P2VO4}. We can only have $r=s=1$. For $1\le j\le 4$, if we partition $S_j$ into blocks as in Example~\ref{eg-G2(3,3)}, we can get $G = 3(3P_2 \vee O_4) \in \mathcal G_4(3,3)$ as follow. The induced vertex labels are 205, 169, 273 respectively.

\begin{figure}[H]
\centerline{\epsfig{file=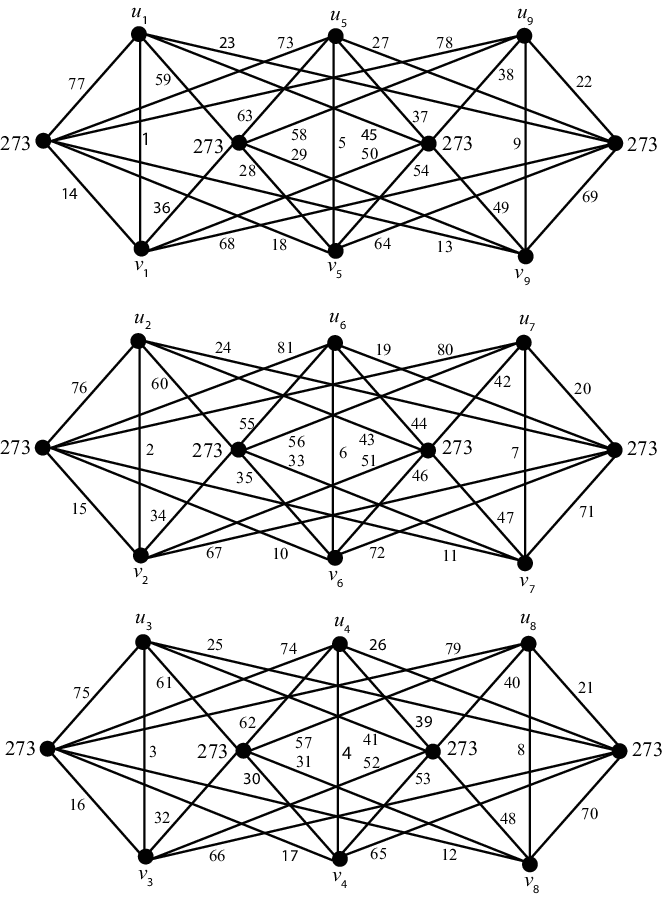, width=8cm}}
\caption{Graph $3(3P_2\vee O_4)\in \mathcal G_4(3,3)$. }\label{fig:3(3P2VO4))}
\end{figure}

\nt If we keep the partitions of $S_1, S_2, S_4$ and partition $S_3$ into blocks $\{f^+(x_{1,3}), f^+(x_{6,3}), f^+(x_{8,3})\}$, $\{f^+(x_{2,3}), f^+(x_{4,3}), f^+(x_{9,3})\}$ and $\{f^+(x_{3,3}), f^+(x_{5,3}), f^+(x_{7,3})\}$, then we get a connected graph as follow with the same induced vertex labels.

\begin{figure}[H]
\centerline{\epsfig{file=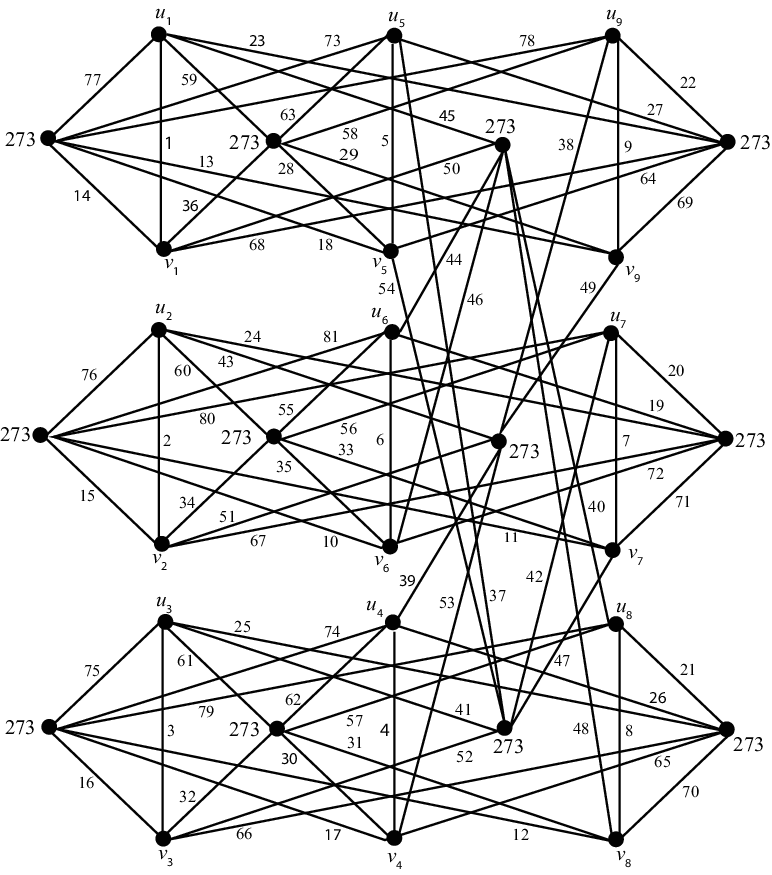, width=9.5cm}}
\caption{A connected graph in $\mathcal G_4(3,3)$.}\label{fig:G4(33)}
\end{figure}

\nt The degree 6 vertices are indicated with their induced vertex label 273. Note that there are many more ways to do the partitioning of $S_1$ to $S_4$ accordingly to get more non-isomorphic graphs all with local antimagic chromatic number 3. \rsq
\end{example}

\nt Recall that for $n\ge 2$, $1\le i\le 2k+1$ and each $j\in [1,2n]\setminus \{2,2n-1\}$, the bijective edge labeling of $(2k+1)(P_2 \vee O_{2n})$  defined above has the property that $f^+(x_{i,j}) = (4k+3)+n(8k+4) =  f^+(x_{k+1, 2}) = f^+(x_{k+1, 2n-1})$. Thus, each $G\in \mathcal G_{2n}(2r+1,2s+1)$ has $(2n-2)(2r+1)$ vertices of degree $4s+2$ each of with is incident to $2s+1$ pairs of edges with labels sum $(4k+3)+n(8k+4)$ and also 2 vertices of degree $4s+2$ each of which is incident to a pair of edges with labels sum $(4k+3)+n(8k+4)$.

\ms\nt For $n\ge 2$, $r,s\ge 1$, and a graph $G\in \mathcal G_{2n}(2r+1,2s+1)$,  choose two distinct vertices of degree $4s+2$, say $x$ and $x'$, such that
\begin{enumerate}[(a)]
\item $xu_i, xv_i$ are edges with labels sum $(4k+3)+n(8k+4)$, while $x'u_{i'}, x'v_{i'}$ are also edges with labels sum $(4k+3)+n(8k+4)$, and
\item $x$ and $x'$ do not have common neighbors.
\end{enumerate}

\nt Delete the edges $xu_i$, $xv_i$, $x'u_{i'}$ and $x'v_{i'}$; and add the edges $xu_{i'}$ and $xv_{i'}$ with labels $f(x'u_{i'})$ and $f(x'v_{i'})$ respectively; and add  the edges $x'u_i$ and $x'v_i$ with labels $f(xu_i)$ and $f(xv_i)$ respectively. Repeat this {\it delete-add} process so long as we get another new graph  not isomorphic to $G$. Let $\mathcal H_{2n}(2r+1,2s+1)$ be the set of all the graphs such obtained. Note that the graphs in $\mathcal H_{2n}(2r+1,2s+1)$ may not be connected.

\begin{theorem}\label{thm-H2n(2r+1,2s+1)} For $n\ge 2, r,s\ge 1$, if $H\in \mathcal H_{2n}(2r+1,2s+1)$, then  $\chi_{la}(H) = 3$. \end{theorem}

\begin{proof} By definition and the discussion above, we immediately have each graph $H \in  \mathcal H_{2n}(2r+1,2s+1)$ also admits a a local antimagic 3-coloring with induced vertex labels as for each graphs   $G\in \mathcal G_{2n}(2r+1,2s+1)$. Moreover, $\chi(H) = 3$. This completes the proof. \end{proof}

\begin{example}\label{eg-H4(33)} Using the $3(3P_2\vee O_3)$ in Example~\ref{eg-G4(33)}, we may let $x$ be the vertex adjacent to $u_1$, $v_1$ and $x'$ be the vertex adjacent to $u_7$, $v_7$, satisfying $f(xu_1) + f(xv_1) =  f(x'u_7) + f(x'v_7) = 91$ and $x,x'$ do not have common neighbors. Note that this graph is not a possible graph of $\mathcal G_4(3,3)$ since two of the degree 6 vertices are obtained by merging vertices with induced vertex labels in $\{f^+(x_{7,4}), f^+(x_{5,1}), f^+(x_{9,1})\}$ and $\{f^+(x_{1,1}), f^+(x_{2,4}), f^+(x_{6,4})\}$ which are not subsets of $S_1$ nor $S_4$.

\begin{figure}[H]
\centerline{\epsfig{file=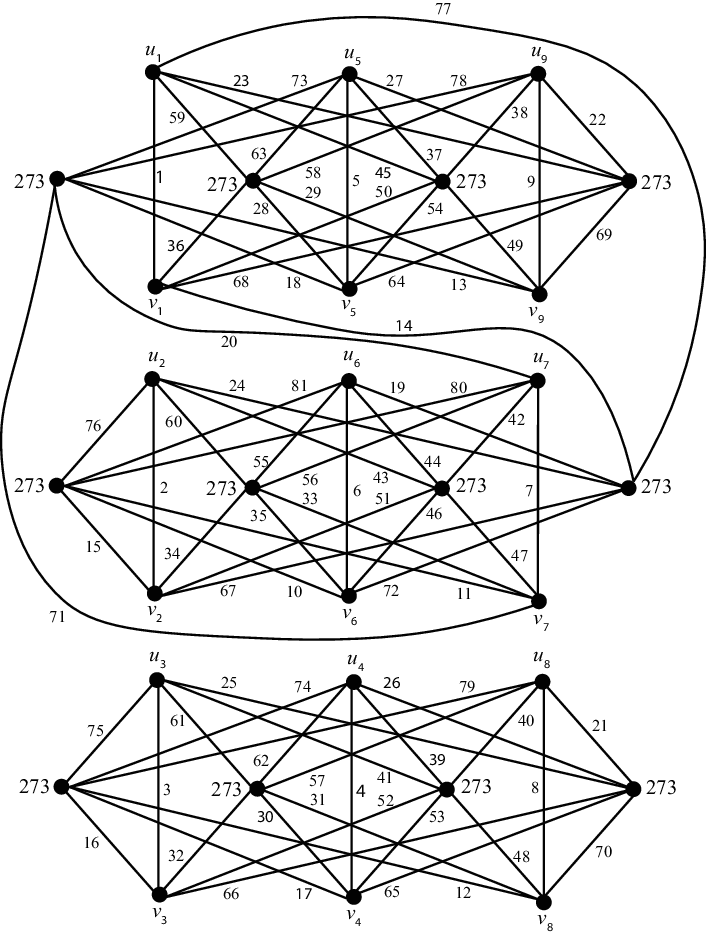, width=8.5cm}}
\caption{A disconnected graph in $\mathcal H_4(3,3)$. }\label{fig:H4(33)1}
\end{figure}

\nt Note that a new connected graph is obtained if we apply the same process to the graph in Figure~\ref{fig:H4(33)1} by choosing the vertices incident to edges with labels $76,15$ and $21,70$. \rsq
\end{example}

\section{$(4n+3)\times (2k+1)$ matrix}\label{subsec:OddP2VOm}

\nt Consider $(2k+1)(P_2\vee O_{2n+1})$ of order $(2k+1)(2n+3)$ and size $(2k+1)(4n+3)$ for $k,n\ge 1$. Thus, for $1\in [1,2k+1]$ and $j\in [1,n]$, we need the following $(4n+3)\times (2k+1)$ matrix with entries in $[1, (4n+3)(2k+1)]$ bijectively. For $n=1$, the required $7\times (2k+1)$ matrix has entries  $f(u_ix_{i,1})$, $f(u_ix_{i,2})$, $f(u_ix_{i,3})$, $f(u_iv_i)$, $f(v_ix_{i,1})$, $f(v_ix_{i,2})$, $f(v_ix_{i,3})$ given by the first 2, middle 3 and last 2 rows.

\[\fontsize{7}{10}\selectfont
\begin{tabu}{|c|[1pt]c|c|c|c|c|c|[1pt]c|}\hline
i & 1 & 2 & 3 & \cdots & 2k & 2k+1 & \mbox{common diff.} \\\tabucline[1pt]{-}
\multirow{2}{1.2cm}{$f(u_ix_{i,1})$}  & 6k+3+ & 6k+2 +  & 6k+1  & \multirow{2}{0.3cm}{$\cdots$} & 4k+4 +& 4k+3 +&  \multirow{2}{0.4cm}{$-1$}  \\
& 2n(4k+2) & 2n(4k+2) & 2n(4k+2) &  & 2n(4k+2) & 2n(4k+2) &  \\\hline
\multirow{2}{1.2cm}{$f(u_ix_{i,2})$}  & 2k+2 + & 2k+3 + & 2k+4 + & \multirow{2}{0.3cm}{$\cdots$} & 4k+1 + & 4k+2 +  &  \multirow{2}{0.4cm}{$+1$} \\
 & 2n(4k+2) & 2n(4k+2) & 2n(4k+2) &  & 2n(4k+2) & 2n(4k+2) &  \\\hline
\vdots & \vdots & \vdots & \vdots & \cdots & \vdots & \vdots  & \vdots \\\hline
\multirow{2}{2cm}{$f(u_ix_{i,2n-2j+1})$} & 6k+3 + & 6k+2 +  & 6k+1  & \multirow{2}{0.3cm}{$\cdots$} & 4k+4 +& 4k+3 +& \multirow{2}{0.4cm}{$-1$} \\
 & (n+j)(4k+2)  &  (n+j)(4k+2) &  (n+j)(4k+2) &   & (n+j)(4k+2) &  (n+j)(4k+2) &  \\\hline
\multirow{2}{2cm}{$f(u_ix_{i,2n-2j+2})$} & 2k+2 + & 2k+3 + & 2k+4 + & \multirow{2}{0.3cm}{$\cdots$} & 4k+1 + & 4k+2 +  &  \multirow{2}{0.4cm}{$+1$} \\
 & (n+j)(4k+2)  &  (n+j)(4k+2) &  (n+j)(4k+2) &   & (n+j)(4k+2) &  (n+j)(4k+2)  &\\\hline
\vdots & \vdots & \vdots & \vdots & \cdots & \vdots & \vdots  & \vdots \\\hline
\multirow{2}{1.6cm}{$f(u_ix_{i,2n+1})$} &  2k+1 + & 2k+ & (2k-1)+  & \multirow{2}{0.3cm}{$\cdots$} & 2 + & 1 + &\multirow{2}{0.4cm}{$-1$}  \\
 &  (n+1)(4k+2) & (n+1)(4k+2) & (n+1)(4k+2) &  & (n+1)(4k+2) & (n+1)(4k+2) & \\\tabucline[1pt]{-}
f(u_iv_i) & 1 & 2 & 3 & \cdots & 2k & 2k+1 & +1 \\\tabucline[1pt]{-}
f(v_ix_{i,1}) & 4k+2 & 4k+1 & 4k & \cdots & 2k+3 & 2k+2 & -1  \\\hline
\vdots & \vdots & \vdots & \vdots & \cdots & \vdots & \vdots  & \vdots  \\\hline
\multirow{2}{1.2cm}{$f(v_ix_{i,2j})$} & 2k+1 + & 2k +  & 2k-1 + & \multirow{2}{0.3cm}{$\cdots$} & 2 + & 1 + &  \multirow{2}{0.4cm}{$-1$}  \\
 & j(4k+2)  &  j(4k+2) &  j(4k+2) &   & j(4k+2) &  j(4k+2) &  \\\hline
\multirow{2}{1.6cm}{$f(v_ix_{i,2j+1})$} & 2k+2 + & 2k+3 +  & 2k+4 + & \multirow{2}{0.3cm}{$\cdots$} & 4k+1 + & 4k +2 +  & \multirow{2}{0.4cm}{$+1$}  \\
 & j(4k+2)  &  j(4k+2) &  j(4k+2) &  & j(4k+2) &  j(4k+2) & \\\hline
\vdots & \vdots & \vdots & \vdots & \cdots & \vdots & \vdots  & \vdots \\\hline
\multirow{2}{1.2cm}{$f(v_ix_{i,2n})$} & 2k+1 + & 2k +  & 2k-1 + & \multirow{2}{0.3cm}{$\cdots$} & 2 + & 1 + &  \multirow{2}{0.4cm}{$-1$}  \\
 & n(4k+2) & n(4k+2) &  n(4k+2) &  &  n(4k+2) & n(4k+2) & \\\hline
\multirow{2}{1.5cm}{$f(v_ix_{i,2n+1})$} & 2k+2 + & 2k+3 +  & 2k+4 + & \multirow{2}{0.3cm}{$\cdots$} & 4k+1 + & 4k +2 +  & \multirow{2}{0.4cm}{$+1$}  \\
 & n(4k+2) & n(4k+2) & n(4k+2) &  & n(4k+2) & n(4k+2) & \\\hline
\end{tabu}
\]

\nt We now have the following observations.
\begin{enumerate}[(1)]
\item For a fixed $j\in [1,n]$,  $\{f(v_i, x_{i, 2j}, f(v_i, x_{i, 2j+1})\;|\; 1\le i\le 2k+1 \} =  [1 + j(4k+2), \ 4k+2+j(4k+2)]. $  Thus, when $j$ runs through $[1,n]$, the integers from $4k+3$ to $(n+1)(4k+2)$ are used.

    Similarly, for a fixed $j\in [1,n]$, $\{f(u_i, x_{i, 2n-2j+1}), f(u_i, x_{i, 2n-2j+2})\;|\;  1\le i\le 2k+1 \}=[2k+2 + (n+j)(4k+2), \ 6k+3 + (n+j)(4k+2)]$.  Thus, when $j$ runs through $[1,n]$, the integers from $2k+2 + (n+1)(4k+2)$ to $2k+1 + (2n+1)(4k+2) = (4n+3)(2k+1)$ are used.

    From rows $f(u_iv_i)$ and $f(v_ix_{i,1})$ with entries 1 to $4k+2$, and rows $f(u_ix_{i,2n+1})$ with entries $1+(n+1)(4k+2)$ to $2k+1+(n+1)(4k+2)$, we see that all integers in $[1, (4n+3)(2k+1)$ are used once.
\item For each column, the sum of the first $2n+2$ entries is $f^+(u_i) = 12n^2k+16nk+6n^2+9n+6k+4$.
\item For each column, the sum of the last $2n+2$ entries is $f^+(v_i) = 4n^2k+8nk + 2n^2 + 5n + 4k+3$.
\item For $S_1 = \{ f^+(x_{i,1})= f(u_ix_{i,1}) + f(v_ix_{i,1})\mid 1\le i\le 2k+1\}$, the elements form an arithmetic sequence with first term $10k+5 + 2n(4k+2)$, last term $6k+5+2n(4k+2)$ and common difference 2. The sum of all these terms is $(2k+1)[(8k+5)+2n(4k+2)]$.
\item For $j\in[2, 2n+1]$, $S_j = \{ f^+(x_{i,j})=f(u_ix_{i,j}) + f(v_ix_{i,j})\mid 1\le i\le 2k+1\} = \{8k+5 + 2n(4k+2), \ldots, 8k+5 + 2n(4k+2)\}$ with multiplicity $2k+1$. The sum of all the terms is $(2k+1)[(8k+5)+2n(4k+2)]$.
\end{enumerate}

\begin{theorem}\label{thm-oddP2VO2n+1} For $n, k\ge 1$, $\chi_{la}((2k+1)P_2 \vee O_{2n+1}) = 3$. \end{theorem}

\begin{proof} Note that $\chi_{la}((2k+1)P_2 \vee O_{2n+1}) \ge \chi((2k+1)P_2 \vee O_{2n}) = 3$. Let $G=(2k+1)(P_2 \vee O_{2n+1})$ for $n, k\ge 1$.  We now define a bijection $f : E(G) \to [1, (4n+3)(2k+1)]$ according to the table above. Now, for each  $j\in [1,2n+1]$, merging the vertices in $\{x_{i,j} \mid 1\le i\le 2k+1\}$, to form new vertex $x_j$ of degree $4k+2$, we get the graph $(2k+1)P_2\vee O_{2n+1}$. From Observations (2) to (5) above, we get that $(2k+1)P_2 \vee O_{2n+1}$ admits a bijective edge labeling $f$ with
\begin{enumerate}[(a)]
\item $f^+(x_j) = (2k+1)[(8k+5)+2n(4k+2)]$,
\item $f^+(u_i) =12n^2k+16nk+6n^2+9n+6k+4$,  and
\item $f^+(v_i) = 4n^2k+8nk + 2n^2 + 5n + 4k+3$,
\end{enumerate}
for $1\le i\le 2k+1$ and $1 \le j\le 2n+1$.

\nt Clearly, for $1\le i\le 2k+1$, $f^+(u_i)> f^+(v_i)$.
Now,
\begin{align*}
(a) - (b) &= 16k^2n-12kn^2+16k^2-6n^2 +12k-5n+1\\
&=(4kn+4k+2n+3)(4k-3n)+4kn+4n+1\\
&> 0 \mbox{ if $4k\ge 3n$.}
\end{align*}
If $3n = 4k+1$, then $(a) - (b) = -(4kn+4k+2n+3) + 4kn + 4n + 1 = -4k + 2n - 2 = -n-1 < 0$. Otherwise, if $3n \ge 4k+2$, we get that $(a) - (b) \le -2(4kn+4k+2n+3) + 4kn + 4n + 1 = -4kn - 8k - 5 < 0$. Thus, $f^+(x_j) \ne f^+(u_i)$. Similarly,
\begin{align*}
(a) - (c) &= 16k^2 n-4kn^2 +16k^2 +8kn-2n^2 +14k-n+2\\
&= (4kn+2n+1)(4k-n)+16k^2+10k+2\\
&> 0 \mbox{ if $n\le 4k$.}
\end{align*}
Otherwise, if $n\ge 4k + 1$ (equivalent, $-n \le -4k-1$), we get that $(a) - (c) \le -4kn-2n + 16k^2 + 10k + 1\le (-4k-1)(4k+2)+16k^2+10k+1 < 0$. Thus, $f^+(x_j) \ne f^+(v_i)$.

\ms\nt Consequently, $f$ is a local antimagic 3-coloring and $\chi_{la}((2k+1)P_2 \vee O_{2n+1}) \le 3$. This completes the proof.
\end{proof}

\begin{example} Consider $n=2, k=4$. The $11\times 9$ matrix and the $9(P_2 \vee O_5)$ with the defined edge labeling are given below. For each $j\in[1,5]$, merge the vertices in $\{x_{i,j}\mid 1\le i\le 9\}$ to get the vertex $x_j$ of degree 18 gives the $9P_2 \vee O_5$ as required. The induced vertex labels of $u_i, v_i, x_j$ are $390, 165, 981$ respectively.

\[
\begin{tabu}{|c|[1pt]c|c|c|c|c|c|c|c|c|}\hline
i & 1 & 2 & 3 & 4 & 5 & 6 & 7 & 8 & 9   \\\tabucline[1pt]{-}
f(u_ix_{i,1}) & 99 & 98 & 97 & 96 & 95 & 94 & 93 & 92 & 91  \\\hline
f(u_ix_{i,2}) & 82 & 83 & 84 & 85 & 86 & 87 & 88 & 89 & 90  \\\hline
f(u_ix_{i,3}) & 81 & 80 & 79 & 78 & 77 & 76 & 75 & 74 & 73   \\\hline
f(u_ix_{i,4}) & 64 & 65 & 66 & 67 & 68 & 69 & 70 & 71 & 72  \\\hline
f(u_ix_{i,5}) & 63 & 62 & 61 & 60 & 59 & 58 & 57 & 56 & 55 \\\tabucline[1pt]{-}
f(u_iv_i) & 1 & 2 & 3 & 4  & 5 & 6 & 7 & 8 & 9 \\\tabucline[1pt]{-}
f(v_ix_{i,1}) & 18 & 17 & 16 & 15 &  14 & 13 &  12 & 11 &  10   \\\hline
f(v_ix_{i,2}) & 27  & 26 & 25 & 24  & 23 & 22 & 21 & 20 & 19   \\\hline
f(v_ix_{i,3})  & 28 & 29 & 30 & 31  & 32 & 33 & 34 & 35 & 36  \\\hline
f(v_ix_{i,4})  & 45  & 44 & 43 & 42 & 41 & 40 & 38 & 38 & 37  \\\hline
f(v_ix_{i,5}) & 46 & 47 & 48 & 49  & 50 & 51 & 52 & 53 & 54 \\\hline
\end{tabu}\]

\begin{figure}[H]
\centerline{\epsfig{file=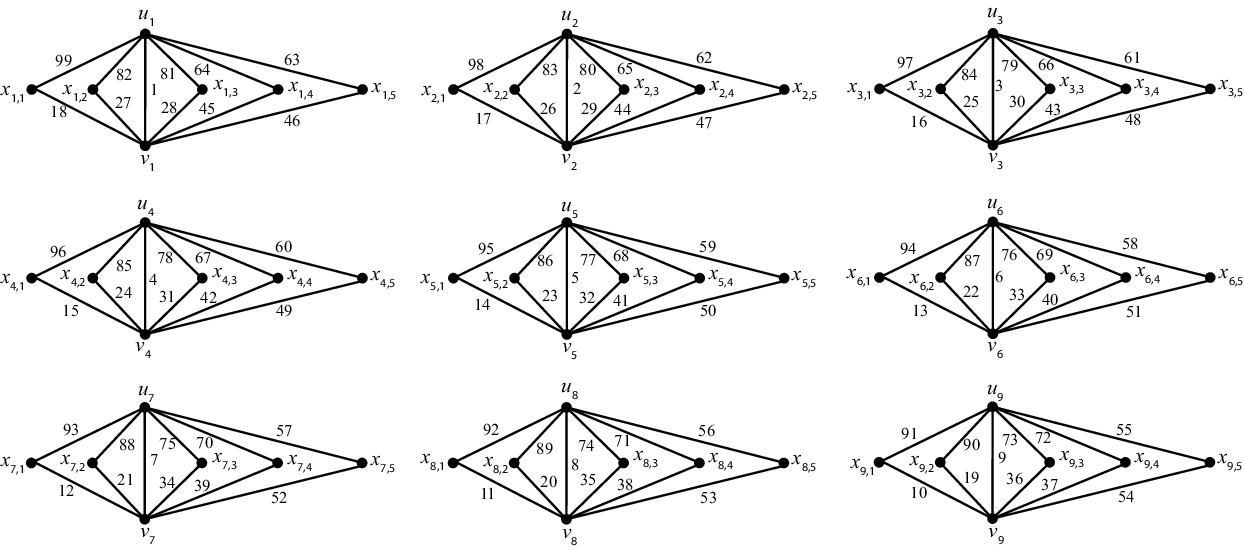, width=15cm}}
\caption{Graph $9(P_2\vee O_5)$.}\label{fig:9(P2VO5))}
\end{figure}\rsq
\end{example}

\nt Similar to Theorem~\ref{thm-G_{2n}(2r+1,2s+1)}, we also define $\mathcal G_{2n+1}(2r+1,2s+1)$ accordingly for $n,r,s\ge 1$.

\begin{theorem}\label{thm-G_{2n+1}(2r+1,2s+1)} For $n,r,s\ge 1$ and $G\in \mathcal G_{2n+1}(2r+1,2s+1)$, $\chi_{la}(G) = 3$. \end{theorem}

\begin{proof}  Using the ideas as in the proof of Theorem~\ref{thm-G_{2n}(2r+1,2s+1)}, we can conclude that $G$ admits a bijective edge labeling with each degree $2(2s+1)$ vertex has induced vertex label $(a) =  (2s+1)[(8k+5)+2n(4k+2)]$ whereas for $1\le i\le 2k+1$, the vertices $u_i$, $v_i$ of degree $2n+2$ have induced vertex labels $(b) = f^+(u_i) = 12n^2k+16nk+6n^2+9n+6k+4$ and $(c) = f^+(v_i) = 4n^2k+8nk + 2n^2 + 5n + 4k+3$. Clearly, $(b) > (c)$. Now,

\begin{align*}
(a) - (b) & = 16kns-12kn^2 -8kn+16ks-6n^2+8ns+2k-5n+10s+1\\
& = (4kn+4k+2n+2)(4s-3n)+4kn+2k+2s+n+1 \\
& > 0 \mbox{ if $4s\ge 3n$.}
\end{align*}
Otherwise, if $4s\le 3n -1$, $(a) - (b) \le 2s-2k-n-1<0$. Thus,  $(a)\ne (b)$. Similarly,
\begin{align*}
(a) - (c) &= 16kns-4kn^2 +16ks+8ns-2n^2 +4k-n+10s+2\\
&= (4kn+2n+1)(4s-n)+16ks+6s+4k+2\\
&> 0 \mbox{ if $4s\ge n$.}
\end{align*}
If $4s \le n-1$ (equivalent, $-n\le -4s-1$), $(a) - (c) \le -4kn-2n+16ks + 6s + 4k + 1\le (-4s-1)(4k+2) + 16ks + 6s + 4k + 1 < 0$. Thus, $(a) \ne (c)$.

\ms\nt Consequently, $\chi_{la}(G)\le 3$. Since $\chi_{la}(G) \ge \chi(G) = 3$, the theorem holds.
\end{proof}


\begin{corollary} For $n, r, s\ge 1$, $\chi_{la}((2r+1)[(2s+1)P_2 \vee O_{2n+1}]) = 3$. \end{corollary}

\begin{proof} Partition $S_1$ so that each of the block  corresponsds to a block of the partition of $S_j, 2\le j\le n$, with elements in the same column of the matrix. \end{proof}

\begin{example} Consider $n=2$, $k=4$. We can only have $r=s=1$. For $1\le j\le 5$, partition $S_j$ into blocks $\{f^+(x_{1,j}), f^+(x_{5,j}), f^+(x_{9j})\}$, $\{f^+(x_{2,j}), f^+(x_{6,j}), f^+(x_{7,j})\}$, $\{f^+(x_{3,j}), f^+(x_{4,j}),$ $f^+(x_{8,j})\}$, we get 6-regular $3(3P_2 \vee O_5) \in \mathcal G_5(3,3)$ as required. The induced vertex labels are $390, 165, 327$ respectively as in Figure~\ref{fig:3(3P2VO5))}. The unlabeled vertices have induced vertex label 327. Similar to Figure~\ref{fig:G4(33)} in Example~\ref{eg-G4(33)}, we can also keep the partition of $S_j, 1\le j\le 4$ and partition $S_5$ into blocks $\{f^+(x_{i,5})\;|\; 1\le i\le 3\}$,  $\{f^+(x_{i,5})\;|\; 4\le i\le 6\}$, and $\{f^+(x_{i,5})\mid 7\le i\le 9\}$ to get a connected graph in $\mathcal G_5(3,3)$. \rsq

\begin{figure}[!ht]
\centerline{\epsfig{file=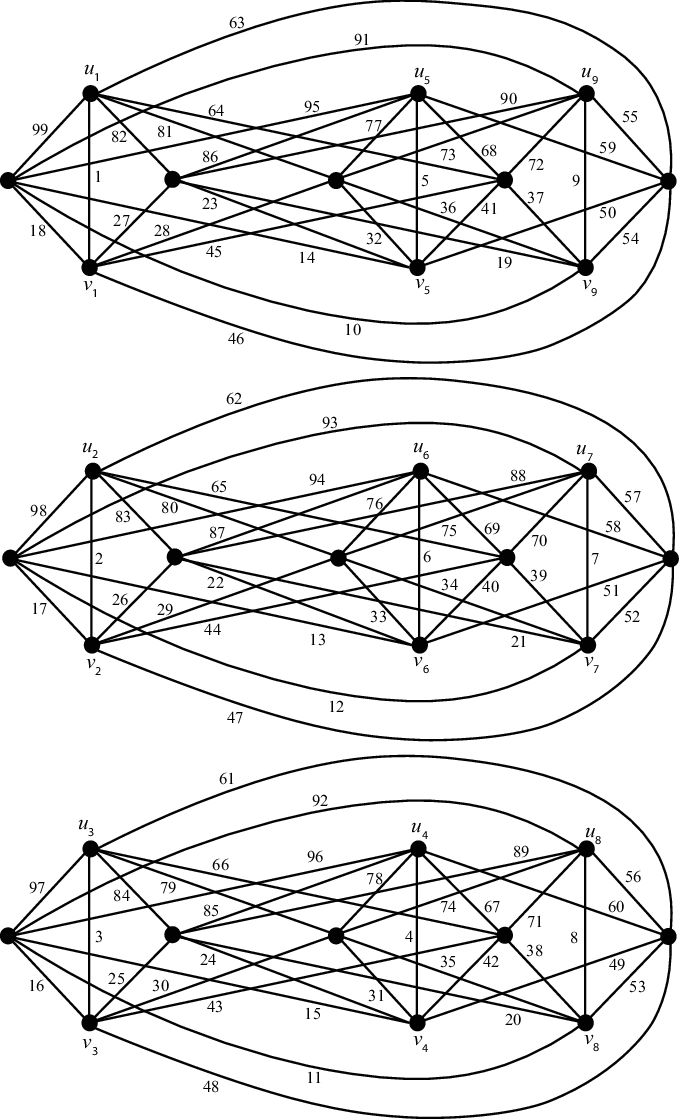, width=8.5cm}}
\caption{Graph $3(3P_2\vee O_5) \in \mathcal G_5(3,3)$.}\label{fig:3(3P2VO5))}
\end{figure}

\end{example}

\nt Similar to Theorem~\ref{thm-H2n(2r+1,2s+1)}, we also have the following theorem without proof.

\begin{theorem}\label{thm-H2n+1(2r+1,2s+1)} For $n,r,s\ge 1$, if $H\in \mathcal H_{2n+1}(2r+1,2s+1)$, then  $\chi_{la}(H) = 3$. \end{theorem}

\nt Similar to Example~\ref{eg-H4(33)}, a (disconnected) graph in $\mathcal H_{2n+1}(2r+1,2s+1)$ can be obtained by applying the  delete-add process to a graph in $\mathcal G_{2n+1}(2r+1,2s+1)$ that has two vertices without common neighbors with a pair of incident edges with equal labels sum. For example, the $3(3P_2 \vee O_5) \in \mathcal G_5(3,3)$ and the two vertices with incident edges labels $54,55$ and $52,57$, respectively.

\begin{corollary}  Each graph (i) $G=(2k+1)P_2\vee O_{4k+1}$, $k\ge 1$, is a connected $(4k+2)$-regular graph of order $8k+3$ and size $(2k+1)(8k+3)$, $G \in \mathcal G_{4s+1}(2r+1,2s+1)$ or $G\in \mathcal H_{4s+1}(2r+1,2s+1)$, $s\ge 1$, is a (possibly disconnected) $(4s+2)$-regular graph of order $(2r+1)(8s+3)$ and size $(2r+1)(2s+1)(8s+3)$  with $\chi_{la}(G) = 3$.   \end{corollary}

\begin{proof} (i) Let $n=2k$, each $G = (2k+1)P_2\vee O_{4k+1}$ is a connected $(4k+2)$-graph. (ii) Let $n = 2s$, each graph in $\mathcal G_{4s+1}(2r+1,2s+1)$ or $\mathcal H_{4s+1}(2r+1,2s+1)$ is a $(4s+2)$ regular graph. \end{proof}

\section{Conclusions and Discussion}

In this paper, we obtained $\chi_{la}((2k+1)P_2 \vee O_m) = 3$ for all $k\ge 1, m\ge 2$. The local antimagic chromatic number of many other families of tripartite graphs are also obtained. As a natural extension, we shall in another paper, make use of matrices of size $(2m+1)\times 2k$ to show that $\chi_{la}((2k)P_2 \vee O_m) = 3$ for all $k\ge 1, m\ge 2$. Consequently, we also obtain the local antimagic chromatic number of many families of bipartite and tripartite graphs.

\end{document}